\documentclass[12pt]{amsart}
\usepackage{amsthm,amsmath,amssymb,amscd,graphics,psfig}

{\end{list}}
{
   \newtheorem{theorem}{Theorem}[section]                     
        
   \newtheorem{lemma}[theorem]{Lemma}

}
{\theoremstyle{definition}

   \newtheorem{remark}[theorem]{Remark}

}
\newcommand{\RR}{{\mathbb{R}}}

\newcommand{\QQ}{{\mathbb{Q}}}

\newcommand{\ZZ}{{\mathbb{Z}}}

\newcommand{\Span}{\operatorname{Span}}

\newcommand{\isom}{\simeq}

\begin{document}
\title{Local Strong Factorization of birational maps} 

\author{Kalle Karu}
\thanks{The author was partially supported by NSF grant DMS-0070678}
\address{Department of Mathematics\\ Harvard University\\ 1 Oxford Street\\
Cambridge, MA 02138\\ USA}
\email{kkaru@math.harvard.edu}

\begin{abstract}
The strong factorization conjecture states that a proper birational
map between smooth algebraic varieties over a field of characteristic
zero can be factored as a sequence of smooth blowups
followed by a sequence of smooth blowdowns. We prove a local version
of the strong factorization conjecture for toric varieties. Combining
this result with the monomialization theorem of S.~D.~Cutkosky, we
obtain a strong factorization theorem for local rings dominated by a
valuation. 
\end{abstract}

\maketitle

\addtocounter{section}{-1}

\section{Introduction}

Let $\phi: X_1\dashrightarrow X_2 $ be a proper birational map between smooth
varieties over a field of characteristic zero. A commutative diagram
\[\begin{array}{rcccl}
 & & Y & & \\
 & \psi_1\swarrow & & \searrow \psi_2 & \\
X_1 & &  \stackrel{\phi}{\dashrightarrow} & & X_2
\end{array}\]
where $\psi_1$ and $\psi_2$ are sequences of blowups of smooth
centers, is called a strong factorization of $\phi$. The existence of
a strong factorization is an open problem in dimension $n=3$ and
higher.

The local version of the strong factorization conjecture replaces the
varieties by local rings dominated by a valuation on their common
fraction field, and the smooth blowups by monoidal transforms along
the valuation. The local strong factorization was proved by
C.~Christensen \cite{Christensen} in dimension $3$ for certain
valuations. A complete proof of the $3$-dimensional case was given by
S.~D.~Cutkosky in \cite{Cutkosky1,Cutkosky2}, where he also made
considerable progress towards proving the conjecture in general. We
prove the local factorization conjecture in any dimension (see
Section~\ref{sec-ring} for notation):

\begin{theorem}\label{thm1} Let $R$ and $S$ be excellent regular local
  rings containing a field $k$ of characteristic zero. Assume that $R$
  and $S$ have a common fraction field $K$ and $\nu$ is a valuation on
  $K$. Then there exists a local ring $T$, obtained from both $R$ and
  $S$ by sequences of monoidal transforms along $\nu$.
\end{theorem}

The toric version of the strong factorization problem considers two
nonsingular fans $\Sigma_1$ and $\Sigma_2$ with the same support and
asks whether there exists a common refinement $\Delta$
\[\begin{array}{rcccl}
 & & \Delta & &  \\
 & \swarrow & & \searrow & \\
\Sigma_1  & & \dashrightarrow & & \Sigma_2
\end{array}\]
obtained from both $\Sigma_1$ and $\Sigma_2$ by sequences of smooth
star subdivisions. Again, this is not known in dimension $3$ or
higher. The local toric version replaces a fan by a single cone and a
smooth star subdivision of the fan by a smooth star subdivision of the
cone together with a choice of one cone in the subdivision. We
assume that the choice is given by a vector $v$ in the cone: we choose
a cone in the subdivision containing $v$. If $v$ has rationally
independent coordinates, then it determines a unique cone in every
subdivision (all cones are rational). We call such a vector $v$ a
{\em valuation} and the subdivision with a choice of a cone a {\em
  subdivision} along the valuation. We prove:

\begin{theorem}\label{thm2} Let $\sigma$ and $\tau$ be nonsingular
  cones, and let $v\in \sigma\cap\tau$ be a vector with rationally
  independent coordinates. Then there exists a nonsingular cone $\rho$
  obtained
  from both $\sigma$ and $\tau$ by sequences of smooth star
  subdivisions along $v$.
\end{theorem}

The proof of Theorem~\ref{thm2} is a generalization of the proof given
by C.~Christensen \cite{Christensen} in dimension $3$. Theorem~\ref{thm1}
follows directly from Theorem~\ref{thm2} and the monomialization
theorem proved by S.~D.~Cutkosky \cite{Cutkosky2}.

Theorem~\ref{thm1} ia also stated in \cite{Cutkosky2}, but the proof
refers to the strong factorization theorem in \cite{Morelli1,
  Abramovich-Matsuki-Rashid} which contains a gap \cite{Matsuki}. We
use the same reduction to the toric case, but replace the reference to
strong factorization by Theorem~\ref{thm2}.

\begin{remark} One can define a more general version of local toric
  factorization. Consider a game between two players $A$ and $B$,
  where the player $A$ subdivides the cone $\tau$ or $\sigma$ and the
  player $B$ chooses one cone in the subdivision (and renames it again
  $\tau$ or $\sigma$). Then the strong factorization conjecture states
  that $A$ always has a winning strategy: after a finite number of
  steps either $\tau=\sigma$ or the interiors of $\tau$ and $\sigma$
  do not intersect. The proof of Theorem~\ref{thm2} given in
  Section~\ref{sec-toric} does not extend to this more general case. A
  positive answer to the global strong factorization conjecture for
  toric varieties would imply that $A$ always has a winning
  strategy. Conversely, a counterexample to the local factorization 
  problem would give a counterexample to the global strong
  factorization conjecture.
\end{remark} 

{\bf Acknowledgments.} I have benefited a great deal from discussions
of the factorization problem with Dan Abramovich, Kenji Matsuki and
Jaros{\l}aw W{\l}odarczyk. It was Jaros{\l}aw's suggestion to look for
a counterexample in dimension $4$ that motivated the current proof.

\section{Local factorization for toric varieties}\label{sec-toric}

Let $N\isom \ZZ^n$ be a lattice and $\sigma$ a rational polyhedral
cone in $N_\RR = N\otimes \RR$ generated by a finite set of vectors
$w_i\in N$
\[ \sigma = \RR_{\geq 0} w_1+\ldots +\RR_{\geq 0} w_m.\]
We say that $\sigma$ is {\em nonsingular} if it can be generated by a
part of a basis of $N$. A nonsingular $m$-dimensional cone has a
unique set of minimal generators $w_1,\ldots, w_m \in N$, and we write 
\[ \sigma = \langle w_1,\ldots,w_m\rangle.\]

We consider nonsingular cones only. When we draw a picture of a cone,
we only show a 
cross-section. Thus, a $3$-dimensional cone is drawn as a triangle.

Let $\sigma = \langle w_1,\ldots,w_n\rangle$ be a nonsingular
$n$-dimensional cone, and let $v\in\sigma$ be a vector $v=c_1
w_1+\ldots + c_n w_n$ such that $c_1,\ldots,c_n$ are linearly
independent over $\QQ$. If $1\leq i< j \leq n$, then precisely one of
the cones  
\[  \langle w_i+w_j, w_1,\ldots, \hat{w}_i,\ldots, w_n\rangle, \qquad 
\langle w_i+w_j, w_1,\ldots, \hat{w}_j,\ldots, w_n\rangle,\]
contains $v$. The cone containing $v$ is called a {\em  star
subdivision of $\sigma$ at $w_i+w_j$ along $v$}. The subdivision is
again a nonsingular cone.  We often denote a
star subdivision of a cone $\sigma$ again $\sigma$, and name its
generators also $w_1,\ldots,w_n$. 

Let us consider the situation of Theorem~\ref{thm2}. It is easy to see
that after star subdividing $\tau$ sufficiently many times we may assume
that $v\in\tau\subset \sigma$. 
We say that a configuration $v\in\tau\subset \sigma$ is {\em
  factorizable} if the statement of 
Theorem~\ref{thm2} holds. We say that $v\in\tau\subset \sigma$ is {\em
  directly factorizable} if the statement of  Theorem~\ref{thm2} holds
with $\rho=\tau$. The vector $v$ is not needed for direct
factorizability.

The following lemma is well-known:

\begin{lemma}\label{lem-dim2} If the dimension  $n=2$, then
  $v\in\tau\subset \sigma$ is   directly factorizable. \qed
\end{lemma}

\begin{lemma}[C.~Christensen \cite{Christensen}]\label{lem-dim3} Let
  $n=3$ and consider $v\in\tau\subset  
  \sigma$, where  
\[ \tau = \langle u_1,u_2,u_3 \rangle, \qquad \sigma= \langle
w_1,w_2,w_3 \rangle\]
are nonsingular cones such that $w_1, u_1, u_2$ are linearly dependent. 
\begin{figure}[h]
\centerline{\psfig{figure=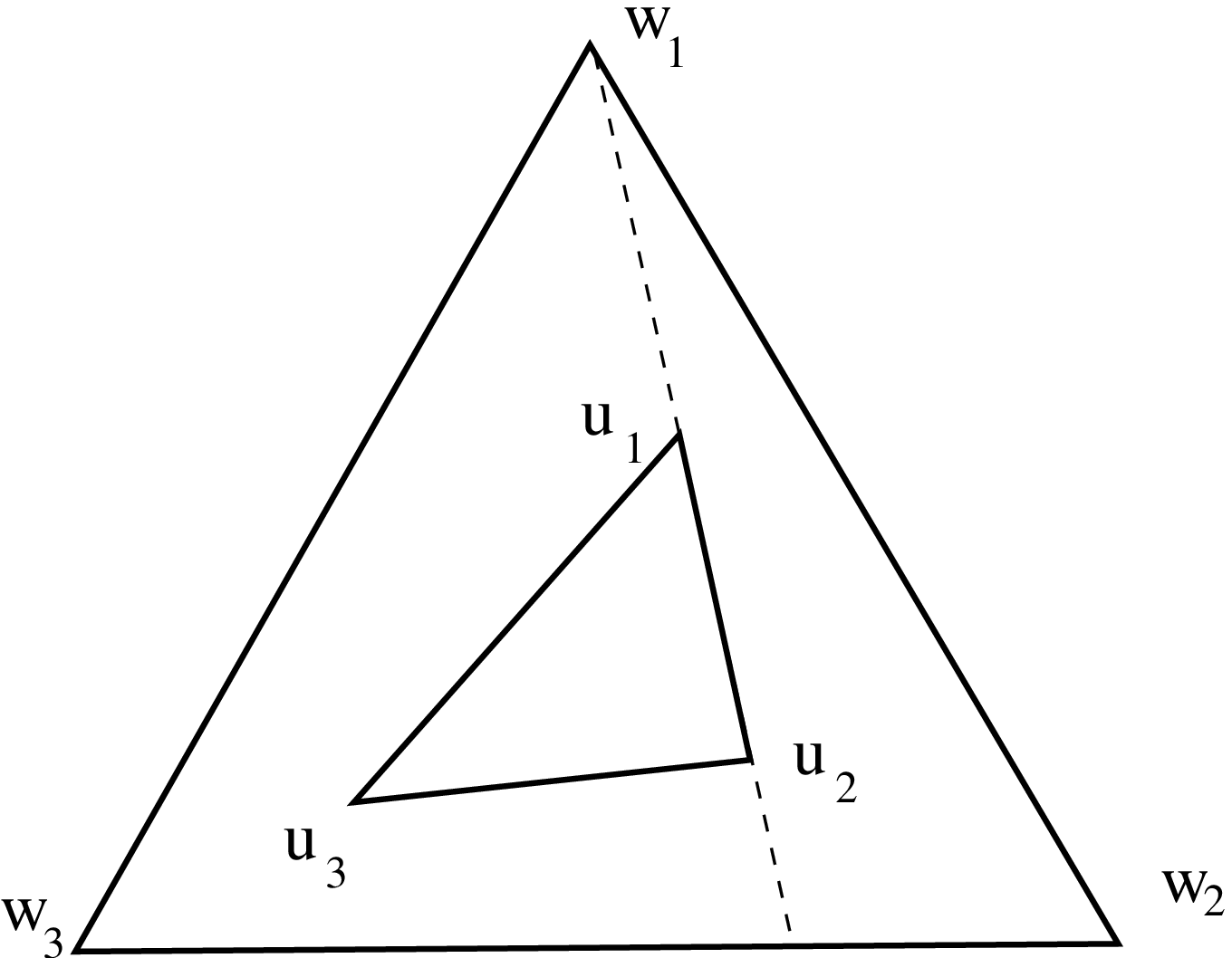,width=2in}}
\end{figure}
Then
$v\in \tau\subset \sigma$ is directly factorizable.
\end{lemma}

{\bf Proof.} Let $\pi: N_\RR \to N_\RR/\RR w_1$ be the quotient
map. We claim that $\pi(\tau)\subset \pi(\sigma)$ are nonsingular
cones with respect to the lattice $\pi(N)$. This is clear for the cone
$\sigma$; for $\tau$ note that the generators $u_1, u_2, u_3$ of $N$
map to generators of $\pi(N)$. (More precisely, $\pi(u_1) = a u'$,
$\pi(u_2) = b u'$, where $u'\in \pi(N)$ is primitive, $\gcd(a,b)=1$.)

Now we apply Lemma~\ref{lem-dim2} to the configuration
$\pi(v)\in \pi(\tau)\subset \pi(\sigma)$. Then after a finite sequence of star
subdivisions of $\sigma$ at vectors lying in $\langle w_2,
w_3\rangle$, we may assume that 
\[ u_3 \in \langle w_1, w_3\rangle, \qquad \langle u_1, u_2\rangle
\subset  \langle w_1, w_2\rangle.\]
\begin{figure}[h]
\centerline{\psfig{figure=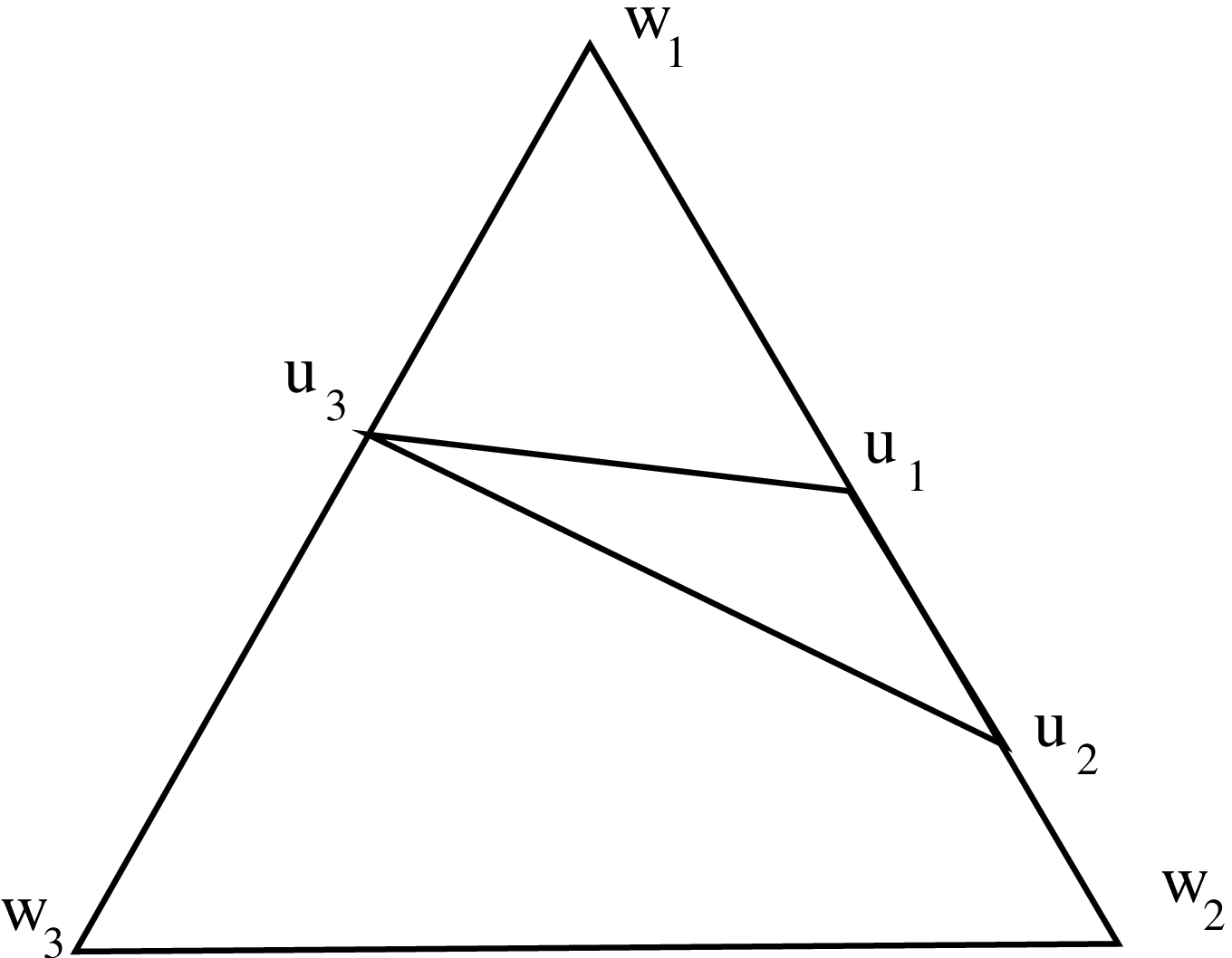,width=2in}}
\end{figure}

If we express $u_3 = \alpha w_1 +\beta w_3$, then it follows from the
nonsingularity of $\tau$ that $\beta=1$. In other words, the cone
$\tau$ lies in the subdivision $\langle w_1+w_3, w_1, w_2\rangle$ of
$\sigma$. Performing a sequence of such star subdivisions, we get to
the situation where $u_3=w_3$. 

\begin{figure}[h]
\centerline{\psfig{figure=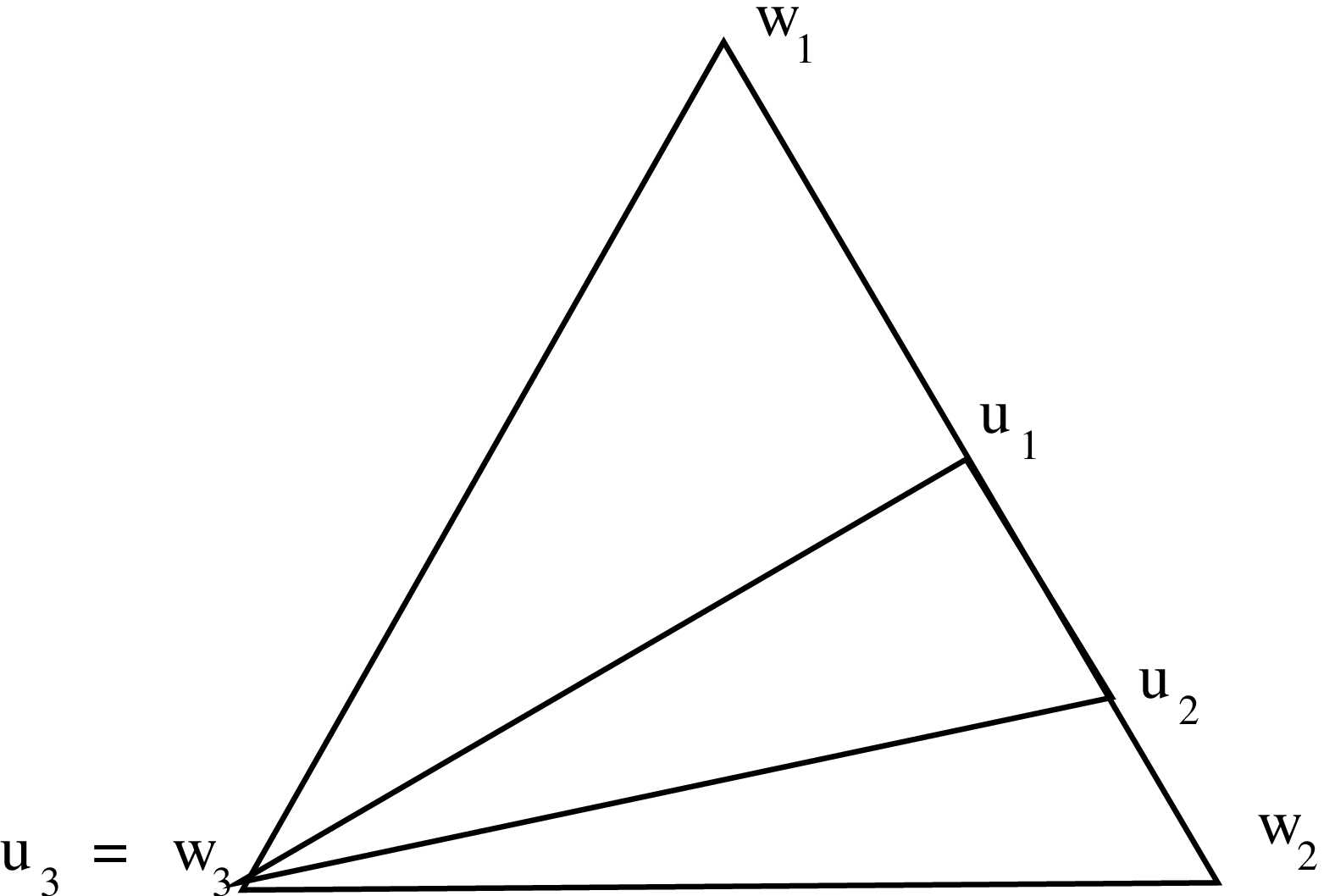,width=2in}}
\end{figure}

Finally,  $\langle u_1, u_2\rangle \subset  \langle w_1, w_2\rangle$
is strongly factorizable by Lemma~\ref{lem-dim2}, thus a
sequence of star subdivisions of $\sigma$ at vectors lying in $\langle
w_1,w_2\rangle$ finishes the proof. \qed

By the previous lemma, to show that $v\in \tau\subset \sigma$ is
factorizable, we have to find a sequence of star subdivisions of
$\tau$ such that the condition of the lemma is satisfied. We prove this
in any dimension.

\begin{lemma}\label{lem-align} Let $n\geq 3$ and consider a
  configuration $v\in   \tau\subset \sigma = \langle w_1,\ldots ,w_n
  \rangle$. There exists 
  a cone $\tau' = \langle u_1,\ldots,u_n \rangle$, obtained from
  $\tau$ by a sequence of smooth star subdivisions along $v$, such
  that $w_1, u_1, u_2$ are linearly dependent. 

Moreover, one can find $\tau'$ such that $w_1, u_1, u_2$ satisfy the
relation
\[ w_1 = u_1-u_2.\]
\end{lemma}

{\bf Proof.} The first part of the proof is again due to C.~Christensen.

Let us start with the case $n=3$ and prove the first half of the
lemma. The algorithm for constructing $\tau'$ is as
follows. Let $\pi: N_\RR \to N_\RR/\RR w_1$ be the projection and let
the generators $u_1, u_2, u_3$ of $\tau$ be ordered so that
$\pi(u_3)\in \pi(\langle u_1, u_2\rangle)$. If $\pi(u_3) \in
\langle\pi(u_1)\rangle$ 
or $\pi(u_3) \in \langle\pi(u_2)\rangle$, then we are done. Otherwise
star subdivide $\tau$ 
at $u_1+u_2$ and repeat.
\begin{figure}[h]
\centerline{\psfig{figure=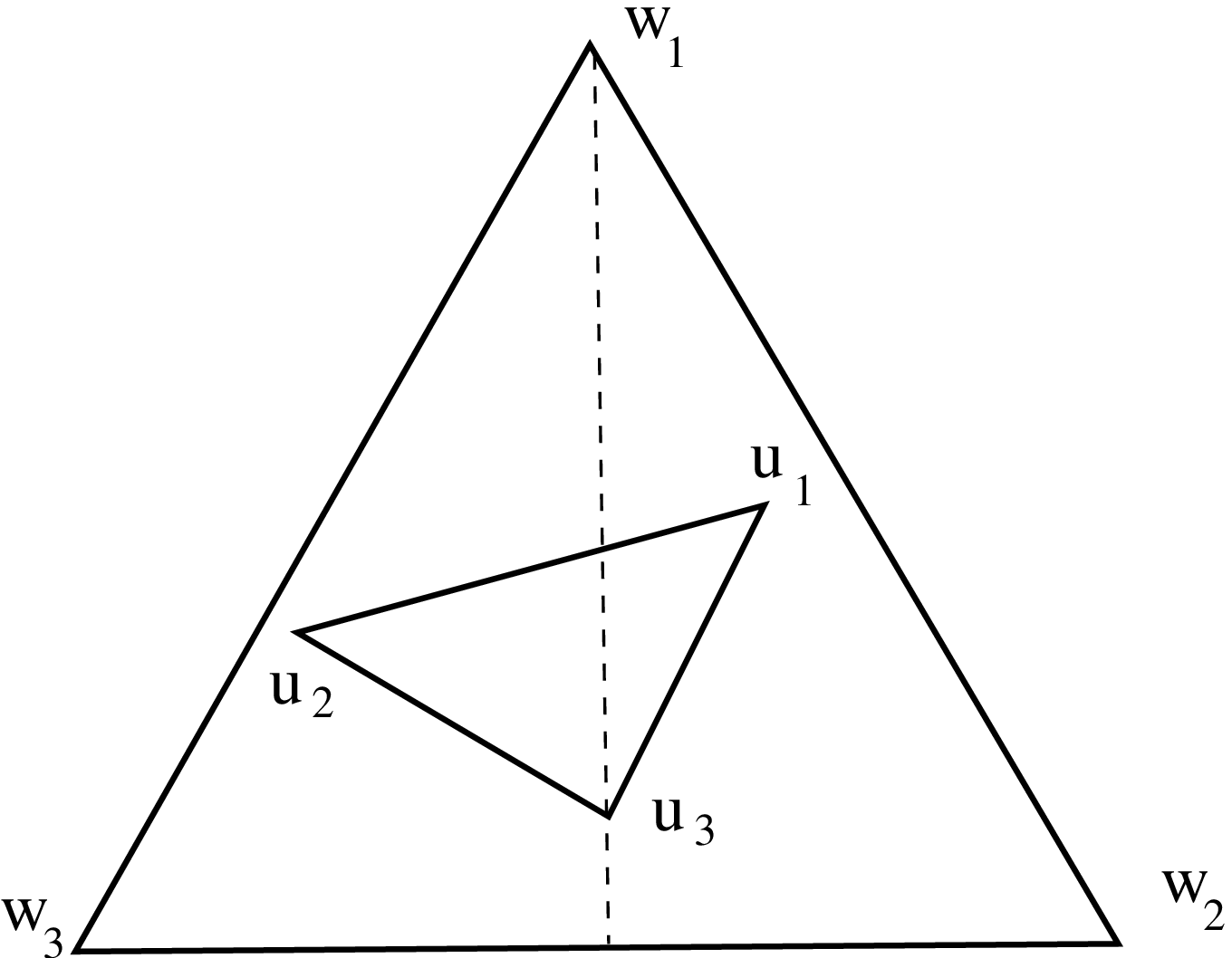,width=2in}}
\end{figure}

To see that this algorithm always terminates, let $a_i, b_i$ be
defined by:
\begin{align*}
w_1 &= a_1 u_1 + a_2 u_2 + a_3 u_3,\\
v &= b_1 u_1 + b_2 u_2 + b_3 u_3.
\end{align*}
Here $a_i\in\ZZ$, $\gcd(a_1,a_2,a_3) = 1$, and $b_i\in\RR$, $b_i>0$.
Then the algorithm can be described as follows. Consider the matrix
\[ \begin{bmatrix} a_1 & a_2 & a_3 \\ b_1 & b_2 & b_3 \end{bmatrix}.\]
If some $a_i = 0$, then we are done. Otherwise, choose columns $i$ and
$j$ such that $a_i$ and $a_j$ have the same sign and subtract the
$i$'th column from the $j$'th if $b_j>b_i$ and $j$'th column from the
$i$'th if $b_i>b_j$. 

Since we always choose columns where $a_i$ and $a_j$ have the same
sign, it is clear that $\max_i |a_i|$ does not increase in this
process, and it suffices to prove that either the algorithm terminates
or $\max_i |a_i|$ drops after a
finite number of steps. Suppose that $|a_3| = \max_i |a_i|$, and
$a_3$ does not change as we run the algorithm. Then $b_3$ also does not
change, and every time we choose columns $i$ and $3$, we subtract
$b_3$ from $b_i$. It is clear that columns $1$ and $2$ can be chosen
only a finite number of times in a row, hence we choose column $3$
infinitely many times. Since we cannot subtract $b_3$ from $b_1$ or
$b_2$ infinitely many times and have a positive result, we get a
contradiction. This proves the first half of the lemma for $n=3$. 

Next let us prove the ``moreover'' part for $n=3$. We start with a
matrix
\[ \begin{bmatrix} a_1 & a_2 & 0 \\ b_1 & b_2 & b_3 \end{bmatrix}.\]

If also $a_2=0$, then by nonsingularity of $\tau$ we have $a_1=1$.
We choose columns $1$ and $2$ the necessary number of times to get
$a_2=-1$: 
\[ \begin{bmatrix} 1 & 0 & 0 \\ b_1 & b_2 & b_3 \end{bmatrix} \to
\begin{bmatrix} 1 & 0 & 0 \\ b_1 - b_2 & b_2 & b_3 \end{bmatrix}
\to\ldots \to \begin{bmatrix} 1 & 0 & 0 \\ b_1 - k b_2 & b_2 & b_3
\end{bmatrix} \to 
\begin{bmatrix} 1 & -1 & 0 \\ b_1 - k b_2 & (k+1)b_2-b_1 & b_3 \end{bmatrix}.
\]

If both $a_1$ and $a_2$ are nonzero then they must have different
signs. Hence, if $\max_i |a_i|=1$ then we are done. Otherwise, since
$\gcd(a_1,a_2)=1$, we may assume that $|a_1|>|a_2|$. We perform star
subdivisions of 
$\tau$ by choosing columns $2$ and $3$ the necessary number of times
to get to the matrix
\[ \begin{bmatrix} a_1 & a_2 & -a_2 \\ b_1 & b_2-k b_3 & (k+1)b_3-b_2
\end{bmatrix}.\] 
After this, we run the algorithm as before. For instance, since $a_1$
and $-a_2$ have the same sign, at the next step we choose columns $1$
and $3$. If we subtract the third column from the first, then  $\max_i
|a_i|$ drops and we are done by induction; otherwise, we subtract the
first column from the third. As before, if $\max_i |a_i|$ does not
decrease, then we are subtracting $b_1$ from $b_2$ or $b_3$ infinitely
many times, and this gives a contradiction.

For $n>3$ we have a matrix
\[ \begin{bmatrix} a_1 & \ldots & a_n \\ b_1 & \ldots & b_n \end{bmatrix}.\]
We can apply the $n=3$ case to the last three columns and achieve
$a_n=0$; then apply the same algorithm to columns $n-3, n-2, n-1$ to
get $a_{n-1} = 0$, and so on, until all but two of the $a_i$ are
nonzero. To prove the second half of the lemma, we apply the $n=3$
case to three columns, including the ones with $a_i\neq 0$. \qed

{\bf Proof of Theorem~\ref{thm2}.} We may assume that $\tau=\langle
u_1,\ldots,u_n\rangle \subset \sigma = \langle w_1,\ldots,w_n\rangle$,
and using Lemma~\ref{lem-align}, we may also assume that $w_1, u_1, u_2$
satisfy the relation
\[ w_1 = u_1-u_2. \qquad (\lozenge)\]
Let $\pi: N_\RR \to N_\RR/\RR w_1$ be the projection. Then $\pi(\tau)
\subset \pi(\sigma)$ are both nonsingular with respect to the lattice
$\pi(N)$. The relation $(\lozenge)$ implies that $\pi(u_1) =
\pi(u_2)$ is a minimal generator of $\pi(\tau)$. In particular, $\pi$
restricts to isomorphisms of cones and lattices:
\begin{gather*}
\pi: \langle u_1, u_3, \ldots,u_n\rangle
\overset{\isom}{\longrightarrow} \pi(\tau), \qquad \pi:
\bigoplus_{i\neq 2} \ZZ u_i \overset{\isom}{\longrightarrow} \pi(N). \\
\pi: \langle u_2, u_3, \ldots,u_n\rangle
\overset{\isom}{\longrightarrow} \pi(\tau), \qquad \pi:
\bigoplus_{i\neq 1} \ZZ u_i \overset{\isom}{\longrightarrow} \pi(N).
\end{gather*}

By induction on the dimension $n$, we have a factorization of $\pi(v)
\in \pi(\tau)\subset \pi(\sigma)$. Unlike the case $n=3$, we may also
have to subdivide $\pi(\tau)$. Consider a star subdivision of
$\pi(\tau)$ at $z = \pi(u_1)+\pi(u_j)$, $j\geq 3$, and define $z'$ and
$z''$ by: 
\[ \pi^{-1}(z)\cap \langle u_2, u_3, \ldots,u_n\rangle = \{z'\},
\qquad \pi^{-1}(z)\cap \langle u_1, u_3, \ldots,u_n\rangle =
\{z''\}.\]

\begin{figure}[h]
\centerline{\psfig{figure=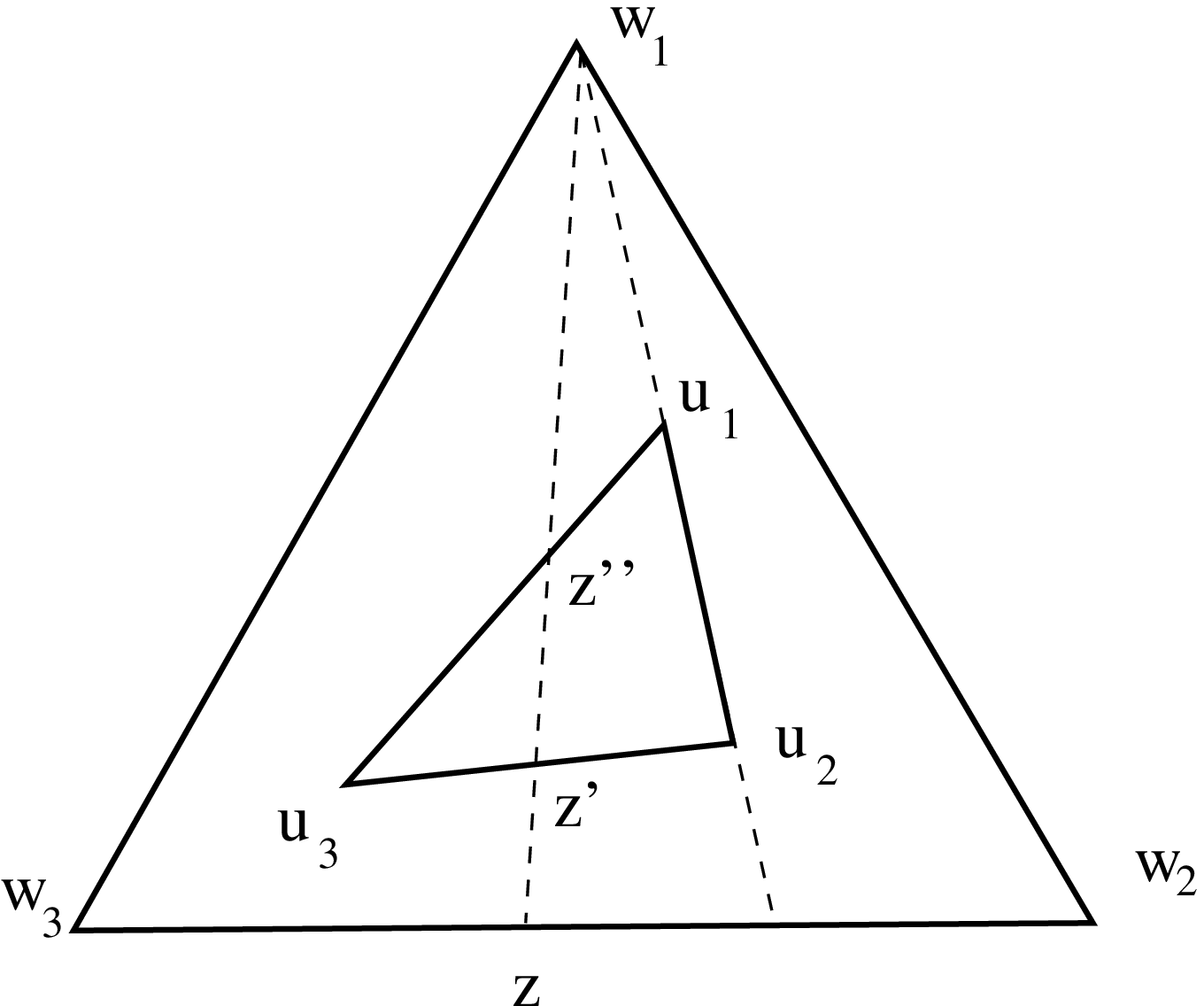,width=2in}}
\end{figure}

Now star subdividing $\tau$ first at $z'$
and then at $z''$ along $v$, the resulting cone again satisfies the relation
$(\lozenge)$ (after possibly reordering the generators), and its
image under $\pi$ is the star subdivision of $\pi(\tau)$ at $z$ along
$\pi(v)$. In other words, every star subdivision of $\pi(\tau)$ can be
lifted to a subdivision of $\tau$. Thus after a finite sequence of
star subdivisions of $\tau$ we may assume that $\pi(v) \in
\pi(\tau)\subset \pi(\sigma)$ is directly factorizable.

The remaining proof is the same as in the $3$-dimensional case. Star
subdividing $\sigma$ at vectors lying in the face $\langle w_2,\ldots, w_n
\rangle$, we may assume that 
\begin{gather*}
 u_i \in  \langle w_1, w_i\rangle, \qquad i=3,\ldots,n, \\
 \langle u_1, u_2\rangle \subset \langle w_1, w_2\rangle.
\end{gather*}
If $u_i = \alpha_i w_1 +\beta_i w_i$, then $\beta_i=1$ for $i\geq 3$,
hence after star subdividing $\sigma$ at vectors lying in the face
$\langle w_1, w_i\rangle$, we may assume that $u_i=w_i$ for $i\geq
3$. Now $\langle u_1, u_2\rangle \subset \langle w_1, w_2\rangle$ are
nonsingular cones, hence directly factorizable by
Lemma~\ref{lem-dim2}. A sequence of star subdivisions of $\sigma$ at
vectors lying in the face $\langle w_1, w_2\rangle$ finishes the proof. \qed

\subsection{The case of $v$ with rationally dependent coordinates.}
J.~W{\l}odarczyk has noted that is makes sense to consider the local
toric factorization problem also for a vector $v$ with rationally
dependent coordinates, and this problem can be reduced to the
rationally independent case. We bring here an argument for such a
reduction. Similar reduction appears in S.~D.~Cutkosky's proof of the
monomialization theorem \cite{Cutkosky2}.

Consider a nonsingular $n$-dimensional cone $\sigma$ and a vector
$v\in\sigma$, with 
possibly rationally dependent coordinates. A star subdivision of $\sigma$
along $v$ is a star subdivision of $\sigma$ and a choice of an
$n$-dimensional cone in the subdivision containing $v$ (i.e., in case
there are more than one such cone, we are free to choose any one of
them). The factorization problem then is: Given two $n$-dimensional
nonsingular cones 
$\sigma, \tau$ and a vector $v\in \sigma\cap\tau$, there exists a
nonsingular cone $\rho$ obtained from both $\tau$ and $\sigma$ by
sequences of star subdivisions along $v$. 
It is clear that the factorization problem has a solution only if the
interiors of $\sigma$ and $\tau$ intersect nontrivially. We assume
that this is the case initially and after every subdivision we
choose a cone containing $v$ such that this condition again
holds. 

An extreme case of the factorization problem is when $v=0$. Then a
factorization along any vector $v'\in\sigma\cap\tau$ (for instance,
$v'$ with rationally independent coordinates) is also a factorization
along $v$. If $v\neq 0$, we reduce the factorization problem to the case of
$v$ with rationally independent coordinates as follows.

The first reduction step is to star subdivide both $\tau$ and $\sigma$
along $v$ to get to the situation where $\tau = \langle u_1,\ldots,
u_n\rangle$, $\sigma = \langle w_1,\ldots,w_n\rangle$ and 
\[ v\in \langle u_1,\ldots, u_m\rangle \cap \langle
w_1,\ldots,w_m\rangle,\]
 such that the coordinates of $v$ with respect to $u_1,\ldots,u_m$
 (hence also with respect to $w_1,\ldots, w_m$) are
 rationally independent. For this write $v=b_1 u_1+\ldots + b_n u_n$
 and consider the vector $(b_1,\ldots,b_n)$ with nonnegative
 entries. It is a simple exercise 
 to show that after a finite sequence of column operations where one
 subtracts $b_i$ from $b_j$ for $b_j\geq b_i$, we get to the vector (after
 reordering the components) $(b_1',\ldots,b_m',0\ldots,0)$, such that
 $b_1',\ldots,b_m'$ are linearly independent over $\QQ$. After a
 similar sequence of star subdivisions of $\sigma$, we get $v = c_1
 w_1+\ldots + c_m w_m$. Note that $\Span(u_1,\ldots,u_m) =
 \Span(w_1,\ldots,w_m)$ is the smallest subspace of $N_\RR$ spanned by
 rational vectors and containing $v$. 

The next step is to use the rationally independent case and factor
$v\in \langle u_1,\ldots, u_m\rangle \cap \langle
w_1,\ldots,w_m\rangle$. Thus after a finite sequence of star
subdivisions of $\sigma$ and $\tau$ we may assume that the two cones
have a common face $\langle u_1,\ldots, u_m\rangle = \langle
w_1,\ldots,w_m\rangle$ containing $v$. After additional subdivisions
of $\tau$ we may also assume that $\tau\subset \sigma$.

The final step is to consider the projection $\pi: N_\RR \to N_\RR/{\RR
w_1}$, and proceed by induction on dimension the same way as in the proof of
Theorem~\ref{thm2}.

\section{Factorization for local rings.} \label{sec-ring}

We recall in this section the monomialization theorem of
S. D. Cutkosky, and then prove Theorem~\ref{thm1}.

Let $(R,\bf{m}_R)$ be a regular local ring of dimension $n$ containing a field
$k$ of characteristic zero, and let $\nu$ be a valuation on the
fraction field of $R$, such that the valuation ring $V$ dominates $R$.
Let $x_1,\ldots,x_m\in R$ be a subset of a system of regular parameters
$x_1,\ldots,x_n$ of $R$. Then the homomorphism 
\[ R \to R'=(R[\frac{x_1}{x_i},\ldots, \frac{x_m}{x_i}])_p, \]
for some $1\leq i\leq m$, and $p$ a prime ideal lying over ${\bf m}_R$, is
called a {\em monoidal transform} of $R$. If $R'$ is again dominated
by the valuation $\nu$, we say that the monoidal transform is a
transform {\em along the valuation $\nu$.} Geometrically, a monoidal
transform is obtained by blowing up a smooth center and localizing at a
point $p$ above ${\bf m}_R$ determined by the valuation. In the following,
we will be interested in monoidal transforms with $m=2$.

Let $R$ and $S$ be two excellent regular local rings of dimension $n$
containing a field $k$ of characteristic zero,
both dominated by a valuation $\nu$ on their common fraction field. 
S.~D.~Cutkosky proved in \cite{Cutkosky1, Cutkosky2} that
after a sequence of monoidal transforms of $R$ and $S$, one can
express a system of regular parameters of $S$ as monomials in regular
parameters of $R$. More precisely, if $\nu$ has rank $1$ and rational
rank $n$ (i.e., the value group can be embedded in $\RR$ and it
contains $n$ rationally independent elements), then after a finite
sequence of monoidal transforms, we may assume that a system
$y_1,\ldots,y_n$ of regular parameters of $S$ can be expressed in
terms of regular parameters $x_1,\ldots, x_n$ of $R$ as:
\begin{gather*}
y_1 = x_1^{a_{1 1}} x_2^{a_{1 2}} \cdots x_n^{a_{1 n}} \\
\ldots     \\
y_n = x_1^{a_{n 1}} x_2^{a_{n 2}} \cdots x_n^{a_{n n}},
\end{gather*}
where $a_{i j}$ are nonnegative integers and $\det(a_{i j}) = \pm
1$. Note that $\nu(y_1),\ldots,\nu(y_n)$ are rationally independent
positive real numbers. If $\nu$ is an arbitrary valuation, then the
matrix $(a_{i j})$ is block diagonal, with $y_i$ corresponding to the
same block having rationally independent values $\nu(y_i)$
(\cite{Cutkosky2}, Theorem~4.4). In the
following proof we will perform monoidal transforms with centers
$(y_i,y_j)$ or $(x_i,x_j)$, with $i$ and $j$ lying in the same block,
hence we may assume that $\nu$ has rank $1$ and rational
rank $n$.

Now let us consider the situation of Theorem~\ref{thm1}. We assume
that an embedding $S \subset R$ is given by monomials as above, and we
have to show that 
after a sequence of monoidal transforms along $\nu$, we get
(renaming parameters) $y_i=x_i$ for $i=1,\ldots,n$. This follows
directly from Theorem~\ref{thm2}, once we express the problem in terms
of cones and subdivisions. 

Let $\sigma =  \langle w_1,\ldots,w_n\rangle \subset N = \ZZ^n$ be a
nonsingular cone, and let $\tau=\langle u_1,\ldots,u_n\rangle \subset
\sigma$ be the cone defined by 
\begin{gather*}
u_1 = a_{1 1} w_1+ a_{2 1} w_2 +\ldots +a_{n 1} w_n\\
\ldots \\
u_n = a_{1 n} w_1+ a_{2 n} w_2 +\ldots +a_{n n} w_n,
\end{gather*}
where $(a_{i j})$ is the matrix of exponents above. Since $\det(a_{i j}) = \pm
1$, the cone $\tau$ is nonsingular. We also let 
\[ v = \nu(y_1) w_1 + \nu(y_2) w_2 +\ldots + \nu(y_n) w_n \]
be a vector $v\in\tau\subset\sigma$. Now one can easily check that
the monoidal transform of $R$ with center $(x_i, x_j)$ along $\nu$
corresponds to the star subdivision of $\tau$ at $u_i+u_j$ along $v$
(which in terms of the matrix $(a_{i j})$ corresponds to adding one
column to another),
and similarly for $S$ and $\sigma$ (in terms of $(a_{i j})$, subtract
one row from another). Applying Theorem~\ref{thm2}, after
a finite sequence of monoidal transforms of $R$ and $S$, the
matrix $a_{i j}$ is the identity matrix, hence $R=S$. \qed

\end{document}